\documentclass[review]{elsarticle}

\usepackage{lineno,hyperref}
\modulolinenumbers[5]
\usepackage{xcolor}
\journal{Journal of \LaTeX\ Templates}
\usepackage{amssymb}

\begin{document}
\newtheorem {Theorem}{\quad Theorem}[section]
\newtheorem{Definition}[Theorem]{\quad Definition}
\newtheorem{Corollary}[Theorem]{\quad Corollary}
\newtheorem{Lemma}[Theorem]{\quad Lemma}
\newtheorem{Proposition}[Theorem]{\quad Proposition}
\newtheorem{Example}[Theorem]{\quad Example}
\newtheorem{Remark}[Theorem]{\quad Remark}
\begin{frontmatter}

\title{Modified different nonlinearities for weakly coupled systems of semilinear effectively damped waves with different time-dependent coefficients in the dissipation terms}

\author{Abdelhamid Mohammed Djaouti}
\ead{djaouti$_{-}$abdelhamid@yahoo.fr; a.mohammeddjaouti@univ-chlef.dz}

\address[mymainaddress]{Faculty of Exact and Computer Sciences, Hassiba Ben Bouali University, Ouled Fares, 02180
	Chlef, Algeria. }

\begin{abstract}
We prove the global existence of small data solution in all space dimension for weakly coupled systems of semi-linear effectively damped wave, with  different time-dependent coefficients in the dissipation terms.  Moreover, nonlinearity terms $ f(t,u) $ and $ g(t,v) $ satisfying some properties of the parabolic equation. We study the problem in several classes of regularity.
\end{abstract}


\end{frontmatter}

\linenumbers

\section{Introduction}\label{intro}
Let us consider the Cauchy problem for semilinear classical damped wave equation with power nonlinearity 
\begin{equation}\label{P1}
u_{tt}-\Delta u+u_{t}=f(u),\,\,\,u(0,x)=u_{0}(x),\,\,\,u_{t}(0,x)=u_{1}(x),
\end{equation} 
where $t\in[0,\infty),\, x\in \mathbb{R}^{n}$ and
\begin{equation}\label{K1}
f(0)=0,\qquad f(u)-f(\tilde{u})\lesssim|u-\tilde{u}|(|u|-|\tilde{u}|)^{p-1}. 
\end{equation}
Having the estimates proved in \cite{Matsumura} for the corresponding homogeneous problem, the authors in \cite{Nakao}  proved for given compactly supported initial data ~$(u_{0},u_{1})\in H^{1}(\mathbb{R}^{n})\times L^{2}(\mathbb{R}^{n})$ and for $p\leq p_{GN}(n):=\frac{n}{n-2}$ if $n\geq 3$ the local (in time) existence of energy solutions $u\in\mathcal{C}([0,T), H^{1}(\mathbb{R}^{n}))\cap \mathcal{C}^{1}([0,T), L^{2}(\mathbb{R}^{n}))$. Moreover, they proved  the global (in time) existence  for small data solutions  by using the technique of potential well and modified potential well. Problem (\ref{P1}) was also devoted  in papers \cite{Fujita,Ikehata0,ikehata,Todorova,Zhang} where Fujita exponent $P_{Fuj}(n):=1+\frac{2}{n}$ has an important role as critical exponent,  which means  that we have global (in time) existence of small data weak solutions for   $p>p_{Fuj}(n),$ while local (in time) existence for $p>1$ with large data.\\
Assuming a time-dependent coefficient in the dissipation term,  we consider first the 
homogeneous problem
\begin{equation}\label{P2}
u_{tt}-\Delta u+b(t)u_{t}=0,\,\,\,
u(0,x)=u_{0}(x),\,\,\,u_{t}(0,x)=u_{1}(x).
\end{equation}
Among other classifications introduced in \cite{Wirth1} and \cite{Wirth2} of the dissipation term $b(t)u_{t}$, we are interested in this paper in the effective case where $ b=b(t) $ satisfies the following properties:
\begin{itemize}
	\item $b$ is a positive and monotonic function with $t b(t)\rightarrow \infty$ as $t\rightarrow \infty$, \label{H2}
	\item $((1+t)^{2}b(t))^{-1} \in L^{1}(0,\infty)$,
	\item  $b\in \mathcal{C}^{3}[0,\infty)$ and
	$	|b^{(k)}(t)| \lesssim \frac{b(t)}{(1+t)^{k}}\,\,\,\mbox{for}\,\,\,k=1,2,3$,
	\item  $\frac{1}{b} \notin L^{1}(0,\infty)$ and there exists a constant $a\in [0,1)$ such that
	$	tb'(t)\leq ab(t)$.
\end{itemize}
Examples of functions belong to this class are: \begin{itemize}
	\item $b(t)=\frac{\mu}{(1+t)^{r}}$ for some  $\mu>0$ and $r\in(-1,1),$
	\item  $b(t)=\frac{\mu}{(1+t)^{r}}(\log(c_{r,\gamma}+t))^{\gamma}$ for some  $\mu>0$ and $\gamma >0$,
	\item  $b(t)=\frac{\mu}{(1+t)^{r}(\log(c_{r,\gamma}+t))^{\gamma}}$ for some  $\mu>0$ and $\gamma >0$.
\end{itemize}	
Here $c_{r,\gamma}$ is a sufficiently large positive constant.\\
In \cite{Reissig1} the authors derived such estimates for solutions to the family of parameter-dependent Cauchy problems
\begin{equation}\label{P3}
u_{tt}-\Delta u+b(t)u_{t}=0,\,\,\,
v(\tau,x)=0,\,\,\,v_{t}(\tau,x)=f(u)(\tau,x).
\end{equation}
Using theses estimates together with Duhamel's Principle the authors proved in the same paper the global existence of small data solutions to the following semilinear Cauchy problem
\begin{equation}\label{P5}
u_{tt}-\Delta u+b(t)u_{t}=f(u),\,\,\,
u(0,x)=u_{0}(x),\,\,\,
u_{t}(0,x)=u_{1}(x),
\end{equation}
where $f(u)$ satisfied condition (\ref{K1}).\\
In 2013, D'Abbicco in \cite{D'Abbicco0} proved the global existence of small data solution for low space dimension  and  derived the decay estimates to the Cauchy problem
$$u_{tt}-\Delta u+b (t)u_{t}=f(t,u), \qquad u(0,x)=u_{0}(x),\,\,\,u_{t}(0,x)=u_{1}(x),$$ 
where $f(0)=0$ and $f(t,v)-f(t,\tilde{v})\lesssim(1+ \int_{0}^{t}\frac{1}{b(r)}dr)^{\gamma}|v-\tilde{v}|(|v|-|\tilde{v}|)^{p-1}.$\\
In this paper we study in all space dimension the   Cauchy problem of weakly coupled system of semilinear effectively damped waves  
\begin{equation}\label{PP1}
\begin{array}{l}
u_{tt}-\Delta u+b_{1}(t)u_{t}=f(t,v), \qquad u(0,x)=u_{0}(x),\,\,\,u_{t}(0,x)=u_{1}(x),\\ v_{tt}-\Delta v+b_{2}(t)v_{t}= g(t,u), \qquad v(0,x)=v_{0}(x),\,\,\,	v_{t}(0,x)=v_{1}(x),
\end{array}
\end{equation}
where
\begin{eqnarray}
\label{K14}&&(1+B_{1}(t,0))^{\frac{1}{\alpha}}\lesssim(1+B_{2}(t,0))\lesssim(1+B_{1}(t,0))^{\beta},\\
\label{K2}&f(0)=0,&\qquad f(t,v)-f(t,\tilde{v})\lesssim(1+B_{1}(t,0))^{\gamma_{1}}|v-\tilde{v}|(|v|-|\tilde{v}|)^{p-1},\\
\label{K15}&g(0)=0,&\qquad g(t,u)-g(t,\tilde{u})\lesssim(1+B_{2}(t,0))^{\gamma_{2}}|u-\tilde{u}|(|u|-|\tilde{u}|)^{q-1},
\end{eqnarray}
for $ B_{1}(t,\tau)=\int_{\tau}^{t}\frac{1}{b_{1}(r)}dr;\,B_{2}(t,\tau)=\int_{\tau}^{t}\frac{1}{b_{2}(r)}dr;\,\alpha,\beta\in\mathbb{R}_{+}^{*}$ and $\gamma_{1},\gamma_{2}\in[-1,\infty).$\\
Recently, K. Nishihara and Y. Wakasugi studied in \cite{Nishihara} the   particular case of (\ref{PP1}), where $b_{1}(t)=b_{2}(t)=1,f(t,v)=|v|^{p}$ and $g(t,u)=|u|^{q}$. Using the weighted energy method they proved the global existence if the inequality
\begin{equation}
\frac{\max\{p;q\}+1}{pq-1}<\frac{n}{2}
\end{equation}
is satisfied. In \cite{Djaouti} and  \cite{Djaouti1} the authors studied the above system with the same nonlinearities assumed in \cite{Nishihara} by taking equivalent coefficients $b_{1}(t)$ and $b_{2}(t)$, or in other word $\alpha=\beta=1$. The global existence for small initial data solutions was proved assuming different classes of regularity of data and for all space dimensions. Considering (\ref{PP1}) in \cite{Djaouti2}, the authors proved a global existence result for a particular case from the set of effective dissipation  terms which is $ b_{1}(t)=\frac{\mu}{(1+t)^{r_{1}}}$  and  $b_{2}(t)=\frac{\mu}{(1+t)^{r_{2}}}$ with the following nonlinearities $ f(t,v)=|v|^{p}$ and $ f(t,u)=|u|^{q}.$ 
\subsection{Notations}
We introduce for $s>0$ and $m \in [1,2)$ the function space
$$\mathcal{A}_{m,s}:=(H^{s}(\mathbb{R}^{n})\cap L^{m}(\mathbb{R}^{n}))\times (H^{s-1}(\mathbb{R}^{n})\cap L^{m}(\mathbb{R}^{n}))$$
with the norm
$$\|(u,v)\|_{\mathcal{A}_{m,s}}:=\|u\|_{H^{s}}+\|u\|_{L^{m}}+ \|v\|_{H^{s-1}}+\|v\|_{L^{m}}.$$
We denote by $\tilde{p}$ and $\tilde{q}$ the modified power nonlinearities of power nonlinearities appeared in (\ref{K2}) and (\ref{K15}). Then
\begin{equation}\label{K5}
\tilde{p}=  \left\{\begin{array}{ccc}
(p-1)\beta+1 & \mbox{if} & \beta\geq1, \\
(p-\frac{m}{2})\beta+\frac{m}{2} & \mbox{if} & 0<\beta<1,
\end{array}
\right.
\end{equation}
and
\begin{equation}\label{K6}
\tilde{q}=  \left\{\begin{array}{ccc}
(q-1)\alpha+1 & \mbox{if} & \alpha\geq1, \\
(q-\frac{m}{2})\alpha+\frac{m}{2} & \mbox{if} & 0<\alpha<1.
\end{array}
\right.
\end{equation}
\begin{Remark}
	If $\max\{\alpha;\beta\}< 1$, then, $(1+B_{1}(t,0)) \approx(1+B_{2}(t,0)).$ This case was studied in previous papers. Then we will restrict ourselves in this work to the remaining cases.
\end{Remark}
\section{Main results}	
We  study the Cauchy problem (\ref{PP1}) in several cases with respect to the regularity of the data in order to cover all space dimonsions, and the modified exponents of power nonlinearities $\,\tilde{p}\,,\,\tilde{q}$ and the parameters $\alpha,\beta,\gamma_{1},\gamma_{2}$. Therefore, we introduce the following classification of regularity:  Data from energy space $s=1$, data from Sobolev spaces with suitable regularity $s\in(1,\frac{n}{2}+1]$ and, finally, large regular data $s>\frac{n}{2}+1.$
\subsection{Data from the energy space}\label{sec:1}  
In this section we are interested in the system (\ref{PP1}), where the data are taken from the function space $\mathcal{A}_{m,1}$. In Theorem \ref{K10} we treat the case where both  modified exponents power nonlinearities $\tilde{p} $ and $\tilde{q} $ are above modified Fujita exponents which are 
$$p_{Fuj,m,\gamma_{1} }:=1+\frac{2m(\gamma_{1} +1)}{n},p_{Fuj,m,\gamma_{2} }:=1+\frac{2m(\gamma_{2} +1)}{n}$$
respectively.
\begin {Theorem}\label{K10}
Let the data  $(u_{0},u_{1}),(v_{0},v_{1})$ are assumed to belong to $\mathcal{A}_{m,1} \times \mathcal{A}_{m,1}$ for $m\in[1,2)$. Moreover, let the modified exponents satisfy
\begin{eqnarray}
\label{K3}\tilde{p}>p_{Fuj,m,\gamma_{1} },\qquad	\label{K4}\tilde{q}>p_{Fuj,m,\gamma_{2} }.
\end{eqnarray}
The exponents $p$ and $q$ of the power nonlinearities satisfy
\begin{equation}\label{Cconditionnenergy}
\begin{array}{ccc}
\frac{2}{m}\leq\min\{p;q\}\leq\max\{p;q\} <\infty  & \mbox{if} & n\leq 2, \\
\frac{2}{m}\leq\min\{p;q\}\leq \max\{p;q\}\leq p_{GN}(n)& \mbox{if} & n>2.
\end{array} 
\end{equation}	
Then, there exists a constant $\epsilon_{0}$ such that if
$$
\|(u_{0},u_{1})\|_{\mathcal{A}_{m,1}}+\|(v_{0},v_{1})\|_{\mathcal{A}_{m,1}}\leq \epsilon_{0},
$$
then there exists a uniquely determined global (in time) energy solution to (\ref{PP1}) in $$\Big(\mathcal{C}\big([0,\infty), H^{1}(\mathbb{R}^{n})\big)\cap\mathcal{C}^1\big([0,\infty), L^{2}(\mathbb{R}^{n})\big)\Big)^2.$$ Furthermore, the solution satisfies the following decay estimate: 
\begin{eqnarray*}
	&&\|\nabla^{j} \partial_{t}^{l}u(t,\cdot)\|_{L^{2}(\mathbb{R}^{n})}\lesssim b_{1}(t)^{-l}\big(1+B_{1}(t,0)\big)^{-\frac{n}{2}\left( \frac{1}{m}-\frac{1}{2}\right)-\frac{j}{2} -l}\\
	&&\qquad\qquad\qquad\qquad\qquad\qquad\qquad\times\big(\|(u_{0},u_{1})\|_{\mathcal{A}_{m,1}}+\|(v_{0},v_{1})\|_{\mathcal{A}_{m,1}}\big),
\end{eqnarray*}
\begin{eqnarray*}
	&&\|\nabla^{j} \partial_{t}^{l}v(t,\cdot)\|_{L^{2}(\mathbb{R}^{n})} \lesssim b_{2}(t)^{-l}\big(1+B_{2}(t,0)\big)^{-\frac{n}{2}\left( \frac{1}{m}-\frac{1}{2}\right)-\frac{j}{2} -l}\\
	&&\qquad\qquad\qquad\qquad\qquad\qquad\qquad\times
	\big(\|(u_{0},u_{1})\|_{\mathcal{A}_{m,1}}+\|(v_{0},v_{1})\|_{\mathcal{A}_{m,1}}\big),
\end{eqnarray*} 
where $j+l=0,1.$
\end{Theorem}
\begin{Remark}
We ramark  for   $\gamma_{1}=\gamma_{2}=0$, that the system (\ref{PP1}) behave like one single equation because the modified  power nonlinearities $\tilde{p}$ and $\tilde{q}$ influenced separatly only by the modified Fujita exponent $p_{Fuj,m}(n)=\frac{2m}{n}+1.$ Then we cannot feel the interplay between the powers of nonlinearities in the existence conditions.
\end{Remark}
\begin{Remark}
The final addmissible ranges for the exponents $p$ and $q$ of power nonlinearities can be fixed using several parameters which are:  $\alpha,\beta,$ the powers $\gamma_1,\gamma_2$, the dimonsion of the space $n$ and  the parameter of additional regularity $m.$ As example for the dimonsion $n=1,$ if we take $0<\beta<1,$ then $\tilde{p}<p.$ We distinguish two cases:
\begin{itemize}
	\item If $ \gamma_1\geq-\frac{1}{2},$ then $p\geq\frac{2}{m}$ is valide for $\tilde{p}>p_{Fuj,m,\gamma_{1} }$ which is equivant to $p>\frac{1}{\beta}\left(2m(\gamma_1+1)-\frac{m}{2} +1\right) +\frac{m}{2}.$
	\item  If $\gamma_1\in[-1,-\frac{1}{2}),$ then the solution existe for $$p>\max\left\lbrace \frac{1}{\beta}\left(2m(\gamma_1+1)-\frac{m}{2}+1 \right)+\frac{m}{2};\frac{2}{m}\right\rbrace .
	$$ 
\end{itemize}  
The general case for the admissible ranges from below can be summarized as follows:\\
$$\begin{tabular}{|l|l|l|}
\hline	 Interplay parameter $ \alpha$ 	&  Nonlinearity parameter $\gamma_1$ & Admissible range for  $p$\\
\qquad & \qquad   &\qquad \\
\hline
$ 0<\beta<1 $ & $ \gamma_1\geq-1+\frac{n}{2} $  &  $ p>\frac{1}{\beta}+\frac{2m(\gamma_1+1)}{n\beta}-\frac{m}{2\beta}+\frac{m}{2} $\\
\cline{2-3}
& $\gamma_1\in[-1,-1+\frac{n}{2})$ & $p>\max\left\lbrace\frac{1}{\beta}+\frac{2m(\gamma_1+1)}{n\beta}-\frac{m}{2\beta}+\frac{m}{2}; \frac{2}{m}\right\rbrace $\\
\hline
$\beta\geq1 $ & $ \gamma_1\geq-1+\frac{n\beta}{2} $  &  $ p> \frac{2m(\gamma_1+1)}{n\beta}+1$\\
\cline{2-3}
& $\gamma_1\in[-1,-1+\frac{n\beta}{2})$ & $ p> \max\left\lbrace \frac{2m(\gamma_1+1)}{n\beta}+1;\frac{2}{m}\right\rbrace $  \\
\hline
\end{tabular}$$
Following similar way one can get the admissible range for $q$ with respect to the parameters $\alpha$ and $\gamma_2$.
\end{Remark}
\begin{Example}\label{K7} Let us choose the dimension space $n=2$, the parameters $\gamma_1=-1,\gamma_2=-\frac{1}{3}$ and the coefficients of the dissipation terms  $b_{1}(t)=(1+t)^{ -\frac{1}{2}}$ and $b_{2}(t)=(1+t)^{ \frac{1}{2}}$
which implies $\beta=\frac{1}{\alpha}=3.$ Using  (\ref{K3}) from previous theorem for $m=2$  we get  $ \tilde{p}>1,\tilde{q}>\frac{7}{3}.$ Theses conditions together with  (\ref{Cconditionnenergy}) after applying (\ref{K5}) and  (\ref{K6}) imply the following admissible range for the exponents of power nonlinearities 
\begin{eqnarray}
p>1,q>\frac{13}{9}.
\end{eqnarray}
\end{Example}
\begin{Example}\label{K8}
If we change  the second coefficient $b_{2}(t)=\frac{(1+t)^{ \frac{1}{2}}}{(\log(e+t))^{\delta}}$, then we consider the Cauchy problem
$$
\begin{array}{l}
u_{tt}-\Delta u+(1+t)^{ -\frac{1}{2}}u_{t}=(1+B_{1}(t,0))^{-1}|v|^{p}, \qquad (u,u_{t})(0,x)=(u_{0},u_{1})(x),\\ v_{tt}-\Delta v+\frac{(1+t)^{ \frac{1}{2}}}{(\log(e+t))^{\delta}}v_{t}= (1+B_{2}(t,0))^{-\frac{1}{3}}|u|^{q}, \qquad (v,v_{t})(0,x)=(v_{0},v_{1})(x),
\end{array}
$$
where $\delta>0$. From the problem we can conclude
$\alpha=1,\beta=3$. In this case we have to garantee $ \tilde{p}>1$ and $\tilde{q}>\frac{7}{3}.$ Finally we conclude the admissible range for the exponents of power nonlinearities  	\begin{eqnarray}
p>1,q>\frac{7}{3}.
\end{eqnarray}
\end{Example}

The case where we have  only one exponent $\tilde{p} $ or $\tilde{q} $ is below modified Fujita exponent, we distinguish four cases with respect to the values of $\alpha$ and $\beta$ as follows:
\begin{enumerate}
\item $\tilde{p}\leq1+\frac{2m(\gamma_{1} +1)}{n},\tilde{q}>1+\frac{2m(\gamma_{2} +1)}{n}$ with  $\min\{\alpha;\beta\}\geq1$ or $\min\{\alpha;\beta\}\leq1\leq\max\{\alpha;\beta\}$.
\item $\tilde{p}>1+\frac{2m(\gamma_{1} +1)}{n},\tilde{q}\leq1+\frac{2m(\gamma_{2} +1)}{n}$ with  $\min\{\alpha;\beta\}\geq1$ or $\min\{\alpha;\beta\}\leq1\leq\max\{\alpha;\beta\}$.
\end{enumerate} 
\begin {Theorem} \label{K11}
Let $m\in[1,2), \alpha\geq 1$ and $\beta>0$. The data  $(u_{0},u_{1}),(v_{0},v_{1})$ are assumed to belong to $\mathcal{A}_{m,1} \times \mathcal{A}_{m,1}$. Moreover, let the modified exponents satisfy
\begin{eqnarray}
\label{K17}\tilde{p}<\frac{2m(\gamma_{1}+1)}{n}+1,\\
\label{K18}\tilde{q}>\frac{2m(\gamma_{2}+1)}{n}+1.
\end{eqnarray}
and
\begin{equation}\label{K9}
\frac{n}{2}>m\left(\frac{\tilde{q}+\alpha+\gamma_{1}\tilde{q}+\gamma_{1}(\alpha-1) +\gamma_{2}}{ \tilde{p}\tilde{q}-1+(\alpha-1)(\tilde{p}-1)}\right).
\end{equation}
The exponents $p$ and $q$ of the power nonlinearities satisfy
\begin{equation}
\begin{array}{ccc}
\frac{2}{m}\leq\min\{p;q\}\leq\max\{p;q\} <\infty  & \mbox{if} & n\leq 2, \\
\frac{2}{m}\leq\min\{p;q\}\leq \max\{p;q\}\leq p_{GN}(n)& \mbox{if} & n>2.
\end{array} 
\end{equation}	
Then, there exists a constant $\epsilon_{0}$ such that if
$$
\|(u_{0},u_{1})\|_{\mathcal{A}_{m,1}}+\|(v_{0},v_{1})\|_{\mathcal{A}_{m,1}}\leq \epsilon_{0},
$$
then there exists a uniquely determined global (in time) energy solution to (\ref{PP1}) in $$\Big(\mathcal{C}\big([0,\infty), H^{1}(\mathbb{R}^{n})\big)\cap\mathcal{C}^1\big([0,\infty), L^{2}(\mathbb{R}^{n})\big)\Big)^2.$$ Furthermore, the solution satisfies the following decay estimates:
\begin{eqnarray*}
&&\|\nabla^{j} \partial_{t}^{l}u(t,\cdot)\|_{L^{2}(\mathbb{R}^{n})} \\
&&\quad\lesssim b_{1}(t)^{-l}\big(1+B_{1}(t,0)\big)^{-\frac{n}{2}\left( \frac{1}{m}-\frac{1}{2}\right)-\frac{j}{2} -l+\kappa(\tilde{p})}
\big(\|(u_{0},u_{1})\|_{\mathcal{A}_{m,1}}+\|(v_{0},v_{1})\|_{\mathcal{A}_{m,1}}\big),
\end{eqnarray*}
\begin{eqnarray*}
&&\|\nabla^{j} \partial_{t}^{l}v(t,\cdot)\|_{L^{2}(\mathbb{R}^{n})}  \\
&&\quad\lesssim b_{2}(t)^{-l}\big(1+B_{2}(t,0)\big)^{-\frac{n}{2}\left( \frac{1}{m}-\frac{1}{2}\right)-\frac{j}{2} -l}
\big(\|(u_{0},u_{1})\|_{\mathcal{A}_{m,1}}+\|(v_{0},v_{1})\|_{\mathcal{A}_{m,1}}\big),
\end{eqnarray*} 
where $j+l=0,1$ and 
\[ \kappa(\tilde{p}) =\gamma_{1}-\frac{n}{2m}(\tilde{p}-1)+1,\] represent the loss of decay in comparison with the corresponding decay estimates for the solution $u$   of the linear Cauchy problem   with vanishing right hand-side.
\end{Theorem} 
\begin{Remark}
If we would choose $\tilde{p}=p_{Fuj,m}(n)$ in condition (\ref{K17}), then we get an  arbitrarily  small loss of decay $\kappa(\tilde{p})=\varepsilon$.
\end{Remark}
We summarize the remaining results for all cases with respect to $\alpha,\beta,\tilde{p}$ and $\tilde{q}$ as follows:
\begin{itemize}
\item 	If we assume in the statement of previous theorem that $\alpha<1$ and $\beta\geq 1,$ the we get instead of (\ref{K9}) the following condition
$$
\frac{n}{2}>m\left(\frac{\tilde{q}+1+\gamma_{1}\tilde{q}+\gamma_{2}+\frac{m}{2}(\alpha-1)(\gamma_{1}+1)}{ \tilde{p}\tilde{q}-1+\frac{m}{2}(\alpha-1)(\tilde{p}-1)}\right).
$$
\item If $\tilde{p}>\frac{2m(\gamma_{1}+1)}{n}+1,\tilde{q}\leq\frac{2m(\gamma_{2}+1)}{n}+1$, then instead of (\ref{K9}) we have to assume
\begin{eqnarray}
\frac{n}{2}&>m\left(\frac{\tilde{p}+\beta+\gamma_2 \tilde{p}+\gamma_{2}(\beta-1) +\gamma_{1}}{ \tilde{p}\tilde{q}-1+(\beta-1)(\tilde{q}-1)}\right) &\quad\mbox{for}\quad  \alpha>0,\quad\beta\geq1,\\
\frac{n}{2}&>m\left(\frac{\tilde{p}+1+\gamma_{2} \tilde{p}+\gamma_{1}+\frac{m}{2}(\beta-1)(\gamma_{2}+1)}{ \tilde{p}\tilde{q}-1+\frac{m}{2}(\beta-1)(\tilde{q}-1)}\right) &\quad\mbox{for}\quad \alpha\geq1,\quad\beta<1.
\end{eqnarray}
\end{itemize}
\subsection{Data from Sobolev spaces with suitable regularity}	\label{Sec2.3}
In this section  the regularity of  data has strong influence on the admissible range of the modified exponents or the exponents of power nonlinearities. For this reason  we assume that the data have different suitable larger regularity, i.e.,
\begin{eqnarray*}
(u_{0},u_{1})\in H^{s_{1}}(\mathbb{R}^{n})\times H^{s_{1}-1}(\mathbb{R}^{n}),&s_{1} \in \Big(1,1+\frac{n}{2}\Big],\\
(v_{0},v_{1})\in H^{s_{2}}(\mathbb{R}^{n})\times H^{s_{2}-1}(\mathbb{R}^{n}),&s_{2} \in \Big(1,1+\frac{n}{2}\Big],
\end{eqnarray*}
with an additional regularity $L^{m}(\mathbb{R}^{n})$, $m \in [1,2)$.
In this section we shall use a generalized (fractional) Gagliardo-Nirenberg inequality used in the papers \cite{hajaieje} and \cite{Kainane}. Furthermore, we shall use  a fractional Leibniz rule and a fractional chain rule which are explained the Appendix.
\begin {Theorem}\label{theoremsamespeeddifferentpowersobolevspacedifferentb}
Let $n\geq4$,  $ s_{1}\in(\max\{1;3+2\gamma_{1}\},\frac{n}{2}+1], s_{2}\in(\max\{1;3+2\gamma_{2}\},\frac{n}{2}+1]$ , $0<s_{2}-s_{1}<1$ and  $\lceil s_{1}\rceil\neq 	\lceil s_{2}\rceil$. The data $(u_{0},u_{1}),(v_{0},v_{1}) $  are supposed to belong to $\mathcal{A}_{m,s_{1}}\times\mathcal{A}_{m,s_{2}}$ with  $m\in[ 1,2)$. Furthermore, we require
 \begin{equation}\label{PP10}
\tilde{p}>\frac{2m}{n}\Big(\frac{s_{1}+1+2\gamma_{1}}{2} \Big) +1,\qquad\tilde{q}>\frac{2m}{n}\Big(\frac{s_{2}+1+2\gamma_{2}}{2} \Big) +1.
\end{equation} 
The exponents $p$ and $q$ of the power nonlinearities satisfy
the conditions
\begin{equation}\label{PP11}
\begin{array}{llcr}
\lceil s_{1}\rceil<p,&\lceil s_{2}\rceil<q &\mbox{if}& n\leq 2s_{1},\\
\lceil s_{1}\rceil<p,&\lceil s_{2}\rceil<q \leq	1+\frac{2}{n-2s_{1}} &\mbox{if}& 2s_{1}<n\leq2s_{2},\\
\lceil s_{1}\rceil<p \leq	1+\frac{2}{n-2s_{2}},& \lceil s_{2}\rceil<q \leq	1+\frac{2}{n-2s_{1}} &\mbox{if}& n>2s_{2}.
\end{array}
\end{equation}
Then, there exists a constant $\epsilon_{0}$ such that if
$$
\|(u_{0},u_{1})\|_{\mathcal{A}_{m,s_1}}+\|(v_{0},v_{1})\|_{\mathcal{A}_{m,s_2}}\leq \epsilon_{0},
$$
then there exists a uniquely determined globally (in time) energy solution to (\ref{PP1}) in  \begin{eqnarray*} && \Big(\mathcal{C}\big([0,\infty), H^{s_{1}}(\mathbb{R}^{n})\big)\cap\mathcal{C}^{1}\big([0,\infty), H^{s_{1}-1}(\mathbb{R}^{n})\big)\Big)  \\ && \qquad \times \Big(\mathcal{C}\big([0,\infty), H^{s_{2}}(\mathbb{R}^{n})\big)\cap\mathcal{C}^{1}\big([0,\infty), H^{s_{2}-1}(\mathbb{R}^{n})\big)\Big).\end{eqnarray*} Furthermore, the solution satisfies for $l=0,1$ the estimates
\begin{eqnarray*}
\||D|^{s_{1}-l}\partial_{t}^{l} u(t,\cdot)\|_{L^{2}(\mathbb{R}^{n})}&\lesssim&b_{1}(t)^{-l} \big(1+B_{1}(t,0)\big)^{-\frac{n}{2}\left( \frac{1}{m}-\frac{1}{2}\right) -l-\frac{s_{1}-l}{2}}\\
&&\qquad\qquad\times\left( \|(u_{0},u_{1})\|_{\mathcal{A}_{m,s_{1}}}+\|(v_{0},v_{1})\|_{\mathcal{A}_{m,s_{2}}}\right),\\		
\||D|^{s_{2}-l}\partial_{t}^{l}v (t,\cdot)\|_{L^{2}(\mathbb{R}^{n})}&\lesssim& b_{2}(t)^{-l}\big(1+B_{2}(t,0)\big)^{-\frac{n}{2}\left( \frac{1}{m}-\frac{1}{2}\right) -l-\frac{s_{2}-l}{2}}\\
&&\qquad\qquad\times \left( \|(u_{0},u_{1})\|_{\mathcal{A}_{m,s_{1}}}+\|(v_{0},v_{1})\|_{\mathcal{A}_{m,s_{2}}}\right).
\end{eqnarray*}
\end{Theorem}
\textbf{Particular cases:}
\begin{itemize}
\item  	If $\beta\geq1$ and $s_{1}\geq3+2\gamma_{1},$ then under the assumptions of Theorem \ref{theoremsamespeeddifferentpowersobolevspacedifferentb}, the condition  $p>\lceil s_{1}\rceil$ implies $\tilde{p}>\frac{2m}{n}\Big(\frac{s_{1}+1+2\gamma_{1}}{2} \Big) +1.$
\item  	If $\alpha\geq1$ and $s_{2}\geq3+2\gamma_{2},$ then under the assumptions of Theorem \ref {theoremsamespeeddifferentpowersobolevspacedifferentb}, the condition  $p>\lceil s_{2}\rceil$ implies $\tilde{q}>\frac{2m}{n}\Big(\frac{s_{2}+1+2\gamma_{2}}{2} \Big) +1.$
\end{itemize}

\subsection{ Large regular data}\label{Sec2.4}
This case has been classified to benefit from the embedding in $ L^{\infty}(\mathbb{R}^{n})$, where the data are supposed to have a high regularity, this means, that
\begin{eqnarray*}
(u_{0},u_{1})\in H^{s_{1}}(\mathbb{R}^{n})\times H^{s_{1}-1}(\mathbb{R}^{n}),&s_{1}> \frac{n}{2}+1,\\
(v_{0},v_{1})\in H^{s_{2}}(\mathbb{R}^{n})\times H^{s_{2}-1}(\mathbb{R}^{n}),&s_{2}>\frac{n}{2}+1.
\end{eqnarray*}
\begin {Theorem} 
\label{differentregularitydifferentb}
Let $n\geq4$, $(u_{0},u_{1}),(v_{0},v_{1})\in \mathcal{A}_{m,s_{1}}\times\mathcal{A}_{m,s_{2}}$, $m\in[ 1,2)$, $\min\{s_{2}; s_{1}\}>\frac{n}{2}+1$, and $s_{1}-s_{2}\in(-1,1)$. Moreover, let
\begin{eqnarray*}
p>s_{1} ,& q>{s}_{2}
\end{eqnarray*}
and 
$$ 
\tilde{p}>\frac{2m}{n}\Big(\frac{s_{1}+1+2\gamma_{1}}{2} \Big) +1,\qquad\tilde{q}>\frac{2m}{n}\Big(\frac{s_{2}+1+2\gamma_{2}}{2} \Big) +1.
$$
Then, there exists a constant $\epsilon_{0}$ such that if
$$
\|(u_{0},u_{1})\|_{\mathcal{A}_{m,s_1}}+\|(v_{0},v_{1})\|_{\mathcal{A}_{m,s_2}}\leq \epsilon_{0},
$$
then there exists a uniquely determined globally (in time) energy solution to (\ref{PP1}) in
\begin{eqnarray*} && \Big(\mathcal{C}\big([0,\infty), H^{s_{1}}(\mathbb{R}^{n})\big)\cap\mathcal{C}^{1}\big([0,\infty), H^{s_{1}-1}(\mathbb{R}^{n})\big)\Big) \\ && \qquad \times \Big(\mathcal{C}\big([0,\infty), H^{{s}_{2}}(\mathbb{R}^{n})\big)\cap\mathcal{C}^{1}\big([0,t], H^{{s}_{2}-1}(\mathbb{R}^{n})\big)\Big).\end{eqnarray*}Furthermore, the solution satisfies for $l=0,1$the   estimates: 	
\begin{eqnarray*}
\||D|^{s_{1}-l}\partial_{t}^{l} u(t,\cdot)\|_{L^{2}(\mathbb{R}^{n})}&\lesssim&b_{1}(t)^{-l} \big(1+B_{1}(t,0)\big)^{-\frac{n}{2}\left( \frac{1}{m}-\frac{1}{2}\right) -l-\frac{s_{1}-l}{2}}\\
&&\qquad\qquad\times\left( \|(u_{0},u_{1})\|_{\mathcal{A}_{m,s_{1}}}+\|(v_{0},v_{1})\|_{\mathcal{A}_{m,s_{2}}}\right),\\		
\||D|^{s_{2}-l}\partial_{t}^{l}v (t,\cdot)\|_{L^{2}(\mathbb{R}^{n})}&\lesssim& b_{2}(t)^{-l}\big(1+B_{2}(t,0)\big)^{-\frac{n}{2}\left( \frac{1}{m}-\frac{1}{2}\right) -l-\frac{s_{2}-l}{2}}\\
&&\qquad\qquad\times \left( \|(u_{0},u_{1})\|_{\mathcal{A}_{m,s_{1}}}+\|(v_{0},v_{1})\|_{\mathcal{A}_{m,s_{2}}}\right).
\end{eqnarray*}
\end{Theorem}
\section{Philosophy of our approach and proofs}
\subsection{Some tools}
First we recall the following result from \cite{Reissig1}.
\begin{Lemma}\label{Bproperties} The primitive  $B(t,\tau)$ satisfies the following properties:
\begin{eqnarray}
\label{P6}B(t,\tau)&\approx& B(t,0) \,\,\,\mbox{for all}\,\,\, \tau\in\Big[0,\frac{t}{2}\Big],\\
\label{P7}	B(\tau,0)&\approx& B(t,0)\,\,\,\mbox{for all}\,\,\, \tau\in\Big[\frac{t}{2},t\Big],
\end{eqnarray}
\begin{equation}\label{P8}
\int_{\frac{t}{2}}^{t}\frac{1}{b(\tau)}\big(1+B(t,\tau)\big)^{-\frac{j}{2}-l}d\tau \lesssim (1+B(t,0))^{1-\frac{j}{2}-l}\log\,(1+B(t,0))^{l}\,\,\,\mbox{for}\,\,\,j+l=0,1.
\end{equation}
\end{Lemma}
In order to use Duhamels principle we need the following results in the proofs of our main results.
\begin{Theorem}\label{linearestimates}
The Sobolev solutions to the Cauchy problem $$
u_{tt}-\Delta u+b(t)u_{t}=0,\,\,\,
u(0,x)=u_{0}(x),\,\,\,u_{t}(0,x)=u_{1}(x)
$$ satisfy the following estimates:\\
For data from the energy space $(s=1)$:
$$\|\nabla^{j} \partial_{t}^{l} u(t,\cdot)\| _{L^{2}}\lesssim (b(t))^{-l}\big(1+B(t,0)\big)^{-\frac{n}{2}\left( \frac{1}{m}-\frac{1}{2}\right)-\frac{j}{2} -l}\|  (u_{0},u_{1})\| _{\mathcal{A}_{m,1}},$$
where $j+l=0,1 1;$\\
for high regular data $(s>1)$:
\begin{eqnarray*}
\| u(t,\cdot)\|_{L^{2}}&\lesssim&  \big(1+B(t,0)\big)^{-\frac{n}{2}\left( \frac{1}{m}-\frac{1}{2}\right) }\| (u_{0},u_{1})\|_{\mathcal{A}_{m,s}},\\
\| u_{t}(t,\cdot)\|_{L^{2}}&\lesssim&  b(t)^{-1}\big(1+B(t,0)\big)^{-\frac{n}{2}\left( \frac{1}{m}-\frac{1}{2}\right) -1}\| (u_{0},u_{1})\|_{\mathcal{A}_{m,s}},\\
\||D|^{s} u(t,\cdot)\|_{L^{2}}&\lesssim&  \big(1+B(t,0)\big)^{-\frac{n}{2}\left( \frac{1}{m}-\frac{1}{2}\right) -\frac{s}{2}}\| (u_{0},u_{1})\|_{\mathcal{A}_{m,s}},\\
\| |D|^{s-1} u_{t}(t,\cdot)\|_{L^{2}}&\lesssim& b(t)^{-1}\big(1+B(t,0)\big)^{-\frac{n}{2}\left( \frac{1}{m}-\frac{1}{2}\right) -\frac{s-1}{2}-1}\| (u_{0},u_{1})\|_{\mathcal{A}_{m,s}}.
\end{eqnarray*}
\end{Theorem}
The proof of this theorem can be concluded from \cite{Wirth1} and \cite{Wirth2}.
\begin{Theorem}\label{nonlinearestimates}
The Sobolev solutions to the parameter-dependent family of Cauchy problems $$
v_{tt}-\Delta v+b(t)v_{t}=0,\,\,\,
v(\tau,x)=0,\,\,\,v_{t}(\tau,x)=v_1(x)
$$ satisfy the following estimates:\\
For	data from the energy space $(s=1)$:
\begin{eqnarray}\label{K12}
\|\nabla^{j}\partial_{t} v(t,\cdot)\| _{L^{2}} \lesssim  b(t)^{-1}b(\tau)^{-l}\big(1+B(t,\tau)\big)^{-\frac{n}{2}\left( \frac{1}{m}-\frac{1}{2}\right) -\frac{j}{2}-l}\| v_1\|_{L^{2}\cap L^{m}},
\end{eqnarray}
where $j+l=0, 1$;\\
for	high regular data $(s>1)$:
\begin{eqnarray}
\|  v(t,\cdot)\| _{L^{2}}&\lesssim & b(\tau)^{-1}\big(1+B(t,\tau)\big)^{-\frac{n}{2}\left( \frac{1}{m}-\frac{1}{2}\right) }\| v_1\|_{H^{s-1}\cap L^{m}},\nonumber\\
\| v_{t}(t,\cdot)\| _{L^{2}}&\lesssim & b(\tau)^{-1}b(t)^{-1}\big(1+B(t,\tau)\big)^{-\frac{n}{2}\left( \frac{1}{m}-\frac{1}{2}\right) -1}\|  v_1\|_{H^{s-1}\cap L^{m}},\nonumber\\
\| |D|^{s} v(t,\cdot)\| _{L^{2}}&\lesssim & b(\tau)^{-1}\big(1+B(t,\tau)\big)^{-\frac{n}{2}\left( \frac{1}{m}-\frac{1}{2}\right) -\frac{s}{2}}\| v_1 \|_{H^{s-1}\cap L^{m}},\nonumber\\
\label{PP4}	\| |D|^{s-1} v_{t}(t,\cdot)\| _{L^{2}}&\lesssim &b(\tau)^{-1}b(t)^{-1}\big(1+B(t,\tau)\big)^{-\frac{n}{2}\left( \frac{1}{m}-\frac{1}{2}\right) -\frac{s-1}{2}-1}\|  v_1\|_{H^{s-1}\cap L^{m}}.
\end{eqnarray}
\end{Theorem}
The proof of this theorem can be concluded from \cite{Reissig1} and \cite{Djaouti3}. 
\subsection{Proofs}
We define the norm of the solution space $X(t)$  by
$$\|(u,v)\|_{X(t)}=\sup_{\tau \in [0,t]}\big\{M_{1}(\tau,u)+M_{2}(\tau,v)\big\},$$
where we shall choose  $M_{1}(\tau,u)$ and $M_{2}(\tau,v)$ with respect to the goals of each theorem.\\
Let $N$ be the mapping on $X(t)$ which is defined  by
$$
N : (u,v) \in X(t) \to N(u,v)=\big(u^{ln}+u^{nl},v^{ln}+v^{nl}\big),
$$
where
\begin{eqnarray*}
u^{ln}(t,x)&:=& E_{1,0}(t,0,x)\ast_{(x)}u_{0}(x)+E_{1,1}(t,0,x)\ast_{(x)}u_{1}(x),\\
u^{nl}(t,x)&:=&\int_{0}^{t}E_{1,1}(t,\tau,x)\ast_{(x)}|v(\tau,x)|^{p}d\tau,\\
v^{ln}(t,x)&:=& E_{2,0}(t,0,x)\ast_{(x)}v_{0}(x)+E_{2,1}(t,0,x)\ast_{(x)}v_{1}(x),\\
v^{nl}(t,x)&:=&\int_{0}^{t}E_{2,1}(t,\tau,x)\ast_{(x)}|u(\tau,x)|^{q}d\tau.
\end{eqnarray*}
We denote by $E_{1,0}=E_{1,0}(t,0,x)$ and $E_{1,1}=E_{1,1}(t,0,x)$ the fundamental solutions to the Cauchy problem
$$u_{tt}-\Delta u+b_{1}(t)u_{t}=0,	\qquad	u(0,x)=u_{0}(x),\qquad	u_{t}(0,x)=u_{1}(x),$$
and by $E_{2,0}=E_{2,0}(t,0,x)$ and $E_{2,1}=E_{2,1}(t,0,x)$ the fundamental solutions to the the Cauchy problem
$$v_{tt}-\Delta v+b_{2}(t)v_{t}=0,	\qquad	v(0,x)=v_{0}(x),\qquad	v_{t}(0,x)=v_{1}(x).$$	
Our aim is to prove the estimates

\begin{equation}\label{PP12}
\begin{array}{l}
\|N(u,v)\|_{X(t)}\\
\qquad \lesssim\|(u_{0},u_{1})\|_{\mathcal{A}_{m,s_1}}+\|(v_{0},v_{1})\|_{\mathcal{A}_{m,s_2}}+\|(u,v)\|_{X(t)}^{p}+\|(u,v)\|_{X(t)}^{q},
\end{array}
\end{equation}
\begin{equation}\label{PP13}
\begin{array}{l}
\|N(u,v)-N(\tilde{u},\tilde{v})\|_{X(t)}\lesssim \|(u,v) - (\tilde{u},\tilde{v})\|_{X(t)}\\
\qquad \times\big(\|(u,v)\|_{X(t)}^{p-1}
+\|(\tilde{u},\tilde{v})\|_{X(t)}^{p-1}
+   \|(u,v)\|_{X(t)}^{q-1}+\|(\tilde{u},\tilde{v})\|_{X(t)}^{q-1}\big).
\end{array}
\end{equation}
We can immediately obtain from	the introduced norm of the solution space $X(t)$ the following inequality:
$$		\|(u^{ln},v^{ln})\|_{X(t)}\lesssim 	 \|(u_{0},u_{1})\|_{\mathcal{A}_{m,s_1}}+\|(v_{0},v_{1})\|_{\mathcal{A}_{m,s_2}}.$$
We complete the proof of all results separately by showing (\ref{PP13}) with the inequality
\begin{equation}\label{P39}
\|(u^{nl},v^{nl})\|_{X(t)}\lesssim \|(u,v)\|_{X(t)}^{p}+\|(u,v)\|_{X(t)}^{q}.
\end{equation} which leads to (\ref{PP12}). \\
\textbf{Proof of Theorem \ref{K10}:} We choose the space of solutions
$$	X(t)=\big( \mathcal{C}([0,t], H^{1})\cap \mathcal{C}^{1}([0,t], L^{2})\big) ^{2},$$
and the following norms  	
\begin{eqnarray*}
M_{1}(\tau,u) &=& (1+B_{1}(\tau,0))^{\frac{n}{2}\left( \frac{1}{m}-\frac{1}{2}\right) }\|u(\tau,\cdot)\|_{L^{2}(\mathbb{R}^{n})}\\
&& +(1+B_{1}(\tau,0))^{\frac{n}{2}\left( \frac{1}{m}-\frac{1}{2}\right)+\frac{1}{2}}\|\nabla u(\tau,\cdot)\|_{L^{2}(\mathbb{R}^{n})}\\
&&+b_{1}(\tau) (1+B_{1}(\tau,0))^{\frac{n}{2}\left( \frac{1}{m}-\frac{1}{2}\right)+1}\| u_{t}(\tau,\cdot)\|_{L^{2}(\mathbb{R}^{n})},\\
M_{2}(\tau,v)&=&  (1+B_{2}(\tau,0))^{\frac{n}{2}\left( \frac{1}{m}-\frac{1}{2}\right) }\|v(\tau,\cdot)\|_{L^{2}(\mathbb{R}^{n})}\\
&& +(1+B_{2}(\tau,0))^{\frac{n}{2}\left( \frac{1}{m}-\frac{1}{2}\right)+\frac{1}{2}}\|\nabla v(\tau,\cdot)\|_{L^{2}(\mathbb{R}^{n})}\\
&&+b_{2}(\tau) (1+B_{2}(\tau,0))^{\frac{n}{2}\left( \frac{1}{m}-\frac{1}{2}\right)+1}\| v_{t}(\tau,\cdot)\|_{L^{2}(\mathbb{R}^{n})}.
\end{eqnarray*}
To prove (\ref{P39}) we need to estimate all  terms appearing in  $\|(u^{nl},v^{nl})\|_{X(t)}$. Let us begin to estimate $\big\|  u^{nl}_{t}(t,\cdot)\big\| _{L^{2}}$. Using (\ref{K12}) with $m=2$ for $\tau\in[\frac{t}{2},t]$  we get
\begin{equation}\label{K13}
\begin{array}{ccc}
\big\| u^{nl}_{t}(t,\cdot)\big\| _{L^{2} } &\lesssim& \displaystyle{\int_{0}^{\frac{t}{2}}}b_{1}(t)^{-1}b_{1}(\tau)^{-1}(1+B_{1}(t,\tau))^{-\frac{n}{2}\left( \frac{1}{m}-\frac{1}{2}\right)-1}\|f(\tau,v)\|_{L^{m} \cap L^{2} }d\tau\\
&&+ \displaystyle{\int_{\frac{t}{2}}^{t}}b_{1}(t)^{-1} b_{1}(\tau)^{-1}(1+B_{1}(t,\tau))^{-1}\| f(\tau,v)\|_{L^{2} }d\tau.
\end{array}
\end{equation}
By a fractional version of Gagliardo-Nirenberg inequality (see proposition \ref{AP1} ) and (\ref{K2}), we obtain
\begin{eqnarray}
\label{K25}&&\| f(\tau,v)\|_{L^{2}} \lesssim (1+B_{1}(\tau,0))^{\gamma_{1}}(1+B_{2}(\tau,0))^{-\frac{n}{2m}p+\frac{n}{4}} 	\|(u,v)\|_{X(t)}^{p},\\
\label{K26}&&\| f(\tau,v)\|_{L^{m}} \lesssim (1+B_{1}(\tau,0))^{\gamma_{1}}(1+B_{2}(\tau,0))^{-\frac{n}{2m}p+\frac{n}{2m}} 	 \|(u,v)\|_{X(t)}^{p},
\end{eqnarray}
where we use condition (\ref{Cconditionnenergy}).
Plugging the last estimates in (\ref{K13}) and using (\ref{K14}), (\ref{P6}) and (\ref{P7}) we get
\begin{eqnarray*}
\big\| u^{nl}_{t}(t,\cdot)\big\| _{L^{2}}&&\lesssim\|(u,v)\|_{X(t)}^{p} \displaystyle{\int_{0}^{\frac{t}{2}}}b_{1}(t)^{-1}b_{1}(\tau)^{-1}(1+B_{1}(t,\tau))^{-\frac{n}{2}\left( \frac{1}{m}-\frac{1}{2}\right)-1} \\&&\qquad\qquad\qquad\qquad\times (1+B_{1}(\tau,0))^{\gamma_{1}}(1+B_{2}(\tau,0))^{-\frac{n}{2m}p+\frac{n}{4}} d\tau\\
&&\qquad\qquad+\|(u,v)\|_{X(t)}^{p} \displaystyle{\int_{\frac{t}{2}}^{t}}b_{1}(t)^{-1} b_{1}(\tau)^{-1}(1+B_{1}(t,\tau))^{-1} \\&&\qquad\qquad\qquad\qquad\times(1+B_{1}(\tau,0))^{\gamma_{1}}(1+B_{2}(\tau,0))^{-\frac{n}{2m}p+\frac{n}{2m}} d\tau\\
&& \lesssim \|(u,v)\|_{X(t)}^{p} \displaystyle{\int_{0}^{\frac{t}{2}}}b_{1}(t)^{-1}b_{1}(\tau)^{-1}(1+B_{1}(t,\tau))^{-\frac{n}{2}\left( \frac{1}{m}-\frac{1}{2}\right)-1}  \\&&\qquad\qquad\qquad\qquad\times(1+B_{1}(\tau,0))^{(-\frac{n}{2m}p+\frac{n}{2m})\beta+\gamma_{1}} d\tau\\
&&\qquad\qquad+\|(u,v)\|_{X(t)}^{p} \displaystyle{\int_{\frac{t}{2}}^{t}}b_{1}(t)^{-1} b_{1}(\tau)^{-1}(1+B_{1}(t,\tau))^{-1} \\&&\qquad\qquad\qquad\qquad\times(1+B_{1}(\tau,0))^{(-\frac{n}{2m}p+\frac{n}{4})\beta+\gamma_{1}} d\tau\\
&& \lesssim \|(u,v)\|_{X(t)}^{p} b_{1}(t)^{-1}(1+B_{1}(t,0))^{-\frac{n}{2}\left( \frac{1}{m}-\frac{1}{2}\right)-1} \\&&\qquad\qquad\qquad\qquad\times\displaystyle{\int_{0}^{\frac{t}{2}}}b_{1}(\tau)^{-1} (1+B_{1}(\tau,0))^{(-\frac{n}{2m}p+\frac{n}{2m})\beta+\gamma_{1}} d\tau\\
&&\qquad\qquad+\|(u,v)\|_{X(t)}^{p}b_{1}(t)^{-1} (1+B_{1}(\tau,0))^{(-\frac{n}{2m}p+\frac{n}{4})\beta+\gamma_{1}}\\&&\qquad\qquad\qquad\qquad\times\displaystyle{\int_{\frac{t}{2}}^{t}} b_{1}(\tau)^{-1}(1+B_{1}(t,\tau))^{-1}  d\tau.
\end{eqnarray*}
We distinguish two cases with respect to the value of $\beta$. If $\beta\geq1$, then we get
\begin{eqnarray*}
\big\| u^{nl}_{t}(t,\cdot)\big\| _{L^{2} } &\lesssim& \|(u,v)\|_{X(t)}^{p} b_{1}(t)^{-1}(1+B_{1}(t,\tau))^{-\frac{n}{2}\left( \frac{1}{m}-\frac{1}{2}\right)-1} \\&&\qquad\qquad\qquad\qquad\times\displaystyle{\int_{0}^{\frac{t}{2}}}b_{1}(\tau)^{-1} (1+B_{1}(\tau,0))^{-\frac{n}{2m}(\tilde{p}-1)+\gamma_{1}} d\tau\\
&&+\|(u,v)\|_{X(t)}^{p}b_{1}(t)^{-1} (1+B_{1}(\tau,0))^{-\frac{n}{2m}(\tilde{p}-1)-\frac{n}{2}\left( \frac{1}{m}-\frac{1}{2}\right) \beta+\gamma_{1}}\\&&\qquad\qquad\qquad\qquad\times\displaystyle{\int_{\frac{t}{2}}^{t}} b_{1}(\tau)^{-1}(1+B_{1}(t,\tau))^{-1}  d\tau\\
&\lesssim& \|(u,v)\|_{X(t)}^{p} b_{1}(t)^{-1}(1+B_{1}(t,\tau))^{-\frac{n}{2}\left( \frac{1}{m}-\frac{1}{2}\right)-1}, 
\end{eqnarray*}
for $\tilde{p}>\frac{2m(\gamma_{1} +1)}{n}+1.$\\
If $0<\beta<1$, then  we get
\begin{eqnarray*}
\big\| u^{nl}_{t}(t,\cdot)\big\| _{L^{2} } &\lesssim& \|(u,v)\|_{X(t)}^{p} b_{1}(t)^{-1}(1+B_{1}(t,\tau))^{-\frac{n}{2}\left( \frac{1}{m}-\frac{1}{2}\right)-1}\\&&\qquad\qquad\qquad\qquad\times \displaystyle{\int_{0}^{\frac{t}{2}}}b_{1}(\tau)^{-1} (1+B_{1}(\tau,0))^{-\frac{n}{2m}(\tilde{p}-1)-\frac{n}{2}\left( \frac{1}{m}-\frac{1}{2}\right)(1-\beta )+\gamma_{1}} d\tau\\
&&+\|(u,v)\|_{X(t)}^{p}b_{1}(t)^{-1} (1+B_{1}(\tau,0))^{-\frac{n}{2m}\tilde{p}+\frac{n}{4}+\gamma_{1}}\\&&\qquad\qquad\qquad\qquad\times\displaystyle{\int_{\frac{t}{2}}^{t}} b_{1}(\tau)^{-1}(1+B_{1}(t,\tau))^{-1}  d\tau\\
&\lesssim& \|(u,v)\|_{X(t)}^{p} b_{1}(t)^{-1}(1+B_{1}(t,\tau))^{-\frac{n}{2}\left( \frac{1}{m}-\frac{1}{2}\right)-1}, 
\end{eqnarray*}
for $\tilde{p}>\frac{2m(\gamma_{1} +1)}{n}+1.$
Finally, we obtain
\begin{equation}\label{K16}
\big\| u^{nl}_{t}(t,\cdot)\big\| _{L^{2} }\lesssim \|(u,v)\|_{X(t)}^{p} b_{1}(t)^{-1}(1+B_{1}(t,\tau))^{-\frac{n}{2}\left( \frac{1}{m}-\frac{1}{2}\right)-1}.
\end{equation}
Analogously, we can prove
\begin{equation}\label{P40}
\big\| \nabla u^{nl}(t,\cdot)\big\| _{L^{2}}\lesssim (1+B_{1}(t,0))^{-\frac{n}{2}\left( \frac{1}{m}-\frac{1}{2}\right) -\frac{1}{2} }\|(u,v)\|_{X(t)}^{p},
\end{equation}
\begin{equation}\label{P41}
\big\| u^{nl}(t,\cdot)\big\| _{L^{2}}\lesssim (1+B_{1}(t,0))^{-\frac{n}{2}\left( \frac{1}{m}-\frac{1}{2}\right) }\|(u,v)\|_{X(t)}^{p}.
\end{equation}
For the second component $v^{nl}$, using Gagliardo-Nirenberg inequality from proposition \ref{AP1} we get
\begin{eqnarray*}
&&\| g(\tau,u)\|_{L^{2}} \lesssim(1+B_{2}(\tau,0))^{\gamma_{1}} (1+B_{1}(\tau,0))^{-\frac{n}{2m}q+\frac{n}{4}} 	\|(u,v)\|_{X(t)}^{q},\\
&&\| g(\tau,u)\|_{L^{m}} \lesssim (1+B_{2}(\tau,0))^{\gamma_{1}}(1+B_{1}(\tau,0))^{-\frac{n}{2m}q+\frac{n}{2m}} 	 \|(u,v)\|_{X(t)}^{q}.
\end{eqnarray*}
Taking account of the last estimates, we can prove similarly  to (\ref{K16})  to (\ref{P41}) the following estimates 
\begin{eqnarray}
\label{P108}
\big\|  v^{nl}_{t}(t,\cdot)\big\| _{L^{2}}&\lesssim&  b_{2}(t)^{-1}	(1+B_{2}(t,0))^{-\frac{n}{2}\left( \frac{1}{m}-\frac{1}{2}\right) -1}\|(u,v)\|_{X(t)}^{q},\\
\label{P42}
\big\| \nabla v^{nl}(t,\cdot)\big\| _{L^{2}}&\lesssim& (1+B_{2}(t,0))^{-\frac{n}{2}\left( \frac{1}{m}-\frac{1}{2}\right) -\frac{1}{2}}\|(u,v)\|_{X(t)}^{q},\\
\label{B19}
\big\| v^{nl}(t,\cdot)\big\| _{L^{2}}&\lesssim& (1+B_{2}(t,0))^{-\frac{n}{2}\left( \frac{1}{m}-\frac{1}{2}\right) }\|(u,v)\|_{X(t)}^{q},
\end{eqnarray}
for $\tilde{q}>\frac{2m(\gamma_{2} +1)}{n}+1.$ 
Finally, (\ref{K16}) to (\ref{B19})  implies (\ref{P39}).\\
The proof of (\ref{PP13}) is completely analogous to the proof of (\ref{PP12}). In this way we complete the proof of Theorem \ref{K10}.\\

\textbf{Proof of Theorem \ref{K11}:} We choose the same space of solutions $X(t)$ and the norm $M_{2}(\tau,v)$ used in the proof of Theorem \ref{K11}.  We modify the norm  $M_{1}(\tau,v)$ as follows:
\begin{eqnarray*}
M_{1}(\tau,u) &=& (1+B_{1}(\tau,0))^{\frac{n}{2}\left( \frac{1}{m}-\frac{1}{2}\right)-\kappa(\tilde{p}) }\|u(\tau,\cdot)\|_{L^{2}(\mathbb{R}^{n})} \\
&&+(1+B_{1}(\tau,0))^{\frac{n}{2}\left( \frac{1}{m}-\frac{1}{2}\right)+\frac{1}{2}-\kappa(\tilde{p})}\|\nabla u(\tau,\cdot)\|_{L^{2}(\mathbb{R}^{n})}\\
&&+b_{1}(\tau) (1+B_{1}(\tau,0))^{\frac{n}{2}\left( \frac{1}{m}-\frac{1}{2}\right)+1-\kappa(\tilde{p})}\| u_{t}(\tau,\cdot)\|_{L^{2}(\mathbb{R}^{n})}.
\end{eqnarray*}
We begin the proof of (\ref{P39}) by the term $\big\| u^{nl}_{t}(t,\cdot)\big\| _{L^{2}}$. Using  (\ref{K12}) with $m=2$ for $\tau\in[\frac{t}{2},t]$  together with Gagliardo-Nirenberg inequality and following the same steps of the proof of (\ref{K16}) we get  
\begin{eqnarray*}
\big\| u^{nl}_{t}(t,\cdot)\big\| _{L^{2} } &\lesssim&  \|(u,v)\|_{X(t)}^{p} b_{1}(t)^{-1}(1+B_{1}(t,\tau))^{-\frac{n}{2}\left( \frac{1}{m}-\frac{1}{2}\right)-1} \\&&\qquad\qquad\qquad\qquad\times\displaystyle{\int_{0}^{\frac{t}{2}}}b_{1}(\tau)^{-1} (1+B_{1}(\tau,0))^{(-\frac{n}{2m}p+\frac{n}{2m})\beta+\gamma_{1}} d\tau\\
&&+\|(u,v)\|_{X(t)}^{p}b_{1}(t)^{-1} (1+B_{1}(\tau,0))^{(-\frac{n}{2m}p+\frac{n}{4})\beta+\gamma_{1}}\\&&\qquad\qquad\qquad\qquad\times\displaystyle{\int_{\frac{t}{2}}^{t}} b_{1}(\tau)^{-1}(1+B_{1}(t,\tau))^{-1} d\tau\\
&\lesssim& \|(u,v)\|_{X(t)}^{p} b_{1}(t)^{-1}(1+B_{1}(t,\tau))^{-\frac{n}{2}\left( \frac{1}{m}-\frac{1}{2}\right)-1+\kappa(\tilde{p})},
\end{eqnarray*}
for $\beta>0.$ Then we have
\begin{equation}\label{K19}
\big\| u^{nl}_{t}(t,\cdot)\big\| _{L^{2} }\lesssim \|(u,v)\|_{X(t)}^{p} b_{1}(t)^{-1}(1+B_{1}(t,\tau))^{-\frac{n}{2}\left( \frac{1}{m}-\frac{1}{2}\right)-1+\kappa(\tilde{p})}.
\end{equation}
In the same way one can prove
\begin{equation}\label{K20}
\big\| \nabla u^{nl}(t,\cdot)\big\| _{L^{2}}\lesssim (1+B_{1}(t,0))^{-\frac{n}{2}\left( \frac{1}{m}-\frac{1}{2}\right) -\frac{1}{2} +\kappa(\tilde{p})}\|(u,v)\|_{X(t)}^{p},
\end{equation}
\begin{equation}\label{K21}
\big\| u^{nl}(t,\cdot)\big\| _{L^{2}}\lesssim (1+B_{1}(t,0))^{-\frac{n}{2}\left( \frac{1}{m}-\frac{1}{2}\right) +\kappa(\tilde{p})}\|(u,v)\|_{X(t)}^{p}.
\end{equation}
Now for $v^{nl}$, we can prove using Gagliardo-Nirenberg inequality and the definition of the solution space $X(t)$ the following estimates 
\begin{eqnarray*}
&&\| g(\tau,u)\|_{L^{2}} \lesssim(1+B_{2}(\tau,0))^{\gamma_{1}} (1+B_{1}(\tau,0))^{-\frac{n}{2m}q+\frac{n}{4}+\kappa(\tilde{p})q} 	\|(u,v)\|_{X(t)}^{q},\\
&&\| g(\tau,u)\|_{L^{m}} \lesssim (1+B_{2}(\tau,0))^{\gamma_{1}}(1+B_{1}(\tau,0))^{-\frac{n}{2m}q+\frac{n}{2m}+\kappa(\tilde{p})q} 	 \|(u,v)\|_{X(t)}^{q}.
\end{eqnarray*} 
Taking account of the last estimates we can prove analogously to (\ref{K19}) to (\ref{K21}) the following
\begin{eqnarray}
\label{K22}
\big\|  v^{nl}_{t}(t,\cdot)\big\| _{L^{2}}&\lesssim&  b(t)^{-1}	(1+B(t,0))^{-\frac{n}{2}\left( \frac{1}{m}-\frac{1}{2}\right) -1}\|(u,v)\|_{X(t)}^{q},\\
\label{K23}
\big\| \nabla v^{nl}(t,\cdot)\big\| _{L^{2}}&\lesssim& (1+B(t,0))^{-\frac{n}{2}\left( \frac{1}{m}-\frac{1}{2}\right) -\frac{1}{2}}\|(u,v)\|_{X(t)}^{q},\\
\label{K24}
\big\| v^{nl}(t,\cdot)\big\| _{L^{2}}&\lesssim& (1+B(t,0))^{-\frac{n}{2}\left( \frac{1}{m}-\frac{1}{2}\right) }\|(u,v)\|_{X(t)}^{q},
\end{eqnarray}
where we use the condition \[ \gamma_{2}-\frac{n}{2m}(\tilde{q}-1) +\kappa(\tilde{p})q\alpha+\varepsilon<-1 \] which is equivalent to condition (\ref{K9}). Consequently, (\ref{K19}) to (\ref{K24})  implies (\ref{P39}) and the proof of Theorem \ref{K11} is completed.  \\
\textbf{Proof of Theorem \ref{theoremsamespeeddifferentpowersobolevspacedifferentb}:} Let us choose the space of solutions
$$ X(t)=\big(\mathcal{C}([0,t], H^{s_{1}})\cap\mathcal{C}^{1}([0,t], H^{s_{1}-1})\big)\times \big(\mathcal{C}([0,t], H^{s_{2}})\cap\mathcal{C}^{1}([0,t], H^{s_{2}-1})\big) $$
with the norm
$$ \|(u,v)\|_{X(t)}=\sup_{\tau \in [0,t]}\big\{M_{1}(\tau,u)+M_{2}(\tau,v)\big\},$$
where
\begin{eqnarray*}
M_{1}(\tau,u) & = & \big(1+B_{1}  (\tau,0)\big)^{\frac{n}{2}\left( \frac{1}{m}-\frac{1}{2}\right) }\|u(\tau,\cdot)\|_{L^{2}(\mathbb{R}^{n})} \\
&&+b_{1}(\tau)\big(1+B_{1}  (\tau,0)\big)^{\frac{n}{2}\left( \frac{1}{m}-\frac{1}{2}\right)+1}\| u_{t}(\tau,\cdot)\|_{L^{2}(\mathbb{R}^{n})}\\
&&+b_{1}(\tau)\big(1+B_{1}  (\tau,0)\big)^{\frac{n}{2}\left( \frac{1}{m}-\frac{1}{2}\right)+\frac{s_{1}-1}{2}+1}\||D|^{s_{1}-1} u_{t}(\tau,\cdot)\|_{L^{2}(\mathbb{R}^{n})} \\
&&+\big(1+B_{1}  (\tau,0)\big)^{\frac{n}{2}\left( \frac{1}{m}-\frac{1}{2}\right)+\frac{s_{1}}{2}}\||D|^{s_{1}} u(\tau,\cdot)\|_{L^{2}(\mathbb{R}^{n})},
\end{eqnarray*}	
and
\begin{eqnarray*}
M_{2}(\tau,v) & = & \big(1+B_{2}  (\tau,0)\big)^{\frac{n}{2}\left( \frac{1}{m}-\frac{1}{2}\right) }\|v(\tau,\cdot)\|_{L^{2}(\mathbb{R}^{n})} \\
&&+b_{2}(\tau)\big(1+B_{2}  (\tau,0)\big)^{\frac{n}{2}\left( \frac{1}{m}-\frac{1}{2}\right)+1}\| v_{t}(\tau,\cdot)\|_{L^{2}(\mathbb{R}^{n})}\\
&&+b_{2}(\tau)\big(1+B_{2}  (\tau,0)\big)^{\frac{n}{2}\left( \frac{1}{m}-\frac{1}{2}\right)+\frac{s_{2}-1}{2}+1}\||D|^{s_{2}-1} v_{t}(\tau,\cdot)\|_{L^{2}(\mathbb{R}^{n})}	\\
&& +\big(1+B_{2}  (\tau,0)\big)^{\frac{n}{2}\left( \frac{1}{m}-\frac{1}{2}\right)+\frac{s_{2}}{2}}\||D|^{s_{2}} v(\tau,\cdot)\|_{L^{2}(\mathbb{R}^{n})}.
\end{eqnarray*}
To prove (\ref{P39}) we  show how to estimate the norms  $\||D|^{s_{1}-1}u^{nl}_{t}(t,\cdot)\|_{L^{2}(\mathbb{R}^{n})}$ and $\||D|^{s_{2}-1}v^{nl}_{t}(t,\cdot)\|_{L^{2}(\mathbb{R}^{n})}$. From the estimate  	(\ref{PP4})  it follows
\begin{eqnarray*}
&&\||D|^{s_{1}-1}u^{nl}_{t}(t,\cdot)\|_{L^{2}(\mathbb{R}^{n})}\\
&&  \lesssim \displaystyle{\int_{0}^{\frac{t}{2}}}b_{1}(\tau)^{-1}b_{1}(t)^{-1}(1+B_{1}(t,\tau))^{-\frac{n}{2}\left( \frac{1}{m}-\frac{1}{2}\right) -\frac{s_{1}-1}{2}-1}\\
&&\qquad\times\| f(\tau,v)\|_{L^{m}(\mathbb{R}^{n})\cap L^{2}(\mathbb{R}^{n})\cap\dot{H}^{s_{1}-1} (\mathbb{R}^{n})}d\tau \\
&&\quad+\displaystyle{\int_{\frac{t}{2}}^{t}}b_{1}(\tau)^{-1}b_{1}(t)^{-1}(1+B_{1}(t,\tau))^{-\frac{n}{2}\left( \frac{1}{m}-\frac{1}{2}\right) -\frac{s_{1}-1}{2}-1}\\
&&\qquad\times\| f(\tau,v)\|_{L^{m}(\mathbb{R}^{n})\cap L^{2}(\mathbb{R}^{n})\cap\dot{H}^{s_{1}-1} (\mathbb{R}^{n})}d\tau.
\end{eqnarray*}
Under the assumptions of Theorem \ref{theoremsamespeeddifferentpowersobolevspacedifferentb} and the choice of the above introduced  norm, for $0\leq\tau\leq t$ the   inequalities (\ref{K25}) and (\ref{K26}) remain true. We calculate the norm \[ \|f(\tau,v)\|_{\dot{H}^{s-1}}.\] 
Using (\ref{ap11}) from the Propositions \ref{ap11} and Proposition \ref{AP1}, we may conclude for $p > \lceil s_{1}-1\rceil$   and $0\leq\tau\leq t$ the following estimate:
\begin{eqnarray*}
\| f(\tau,v)\|_{\dot{H}^{s_{1}-1}} &\lesssim& (1+B_{1}(\tau,0))^{\gamma_{1}} \big\|  v(\tau,\cdot)\big\|^{p-1} _{L^{q_{1}}} \big\| |D|^{s_{1}-1} (\tau,\cdot)\big\| _{L^{q_{2}}}\\
&\lesssim&(1+B_{1}(\tau,0))^{\gamma_{1}}\big\|  v(\tau,\cdot)\big\|^{(p-1)(1-\theta_{1})} _{L^{2}}\\
&&\qquad\times\big\|  |D|^{s_{2}}v(\tau,\cdot)\big\|^{(p-1)\theta_{1}}_{L^{2}}\big\| v(\tau,\cdot)\big\|^{1-\theta_{2}}_{L^{2}}\big\|  |D|^{s_{2}}v(\tau,\cdot)\big\|^{\theta_{2}} _{L^{2}}\\
&\lesssim& (1+B_{1}(\tau,0))^{\gamma_{1}}(1+B_{2}(\tau,0))^{-\frac{n}{2m}p+\frac{n}{4}-\frac{s_{1}-1}{2}}\|(u,v)\|^{p}_{X(t)},
\end{eqnarray*}
where $$\frac{p-1}{q_{1}}+\frac{1}{q_{2}}=\frac{1}{2},\,\,\,\theta_{1}=\frac{n}{s}\Big(  \frac{1}{2}-\frac{1}{q_{1}}\Big) \in [0,1],\,\,\,\theta_{2}=\frac{n}{s_{2}}\Big(  \frac{1}{2}-\frac{1}{q_{2}}\Big) + \frac{s_{1}-1}{s_{2}} \in \Big[  \frac{s_{1}-1}{s_{2}},1\Big].$$
To satisfy the last conditions for the parameters $\theta_{1}$ and $\theta_{2}$ we choose $q_{2}=\frac{2n}{n-2}$ and $q_{1}=n(p-1)$. This choice implies the condition \[  1+\frac{2}{n}  \leq p\leq 1+\frac{2}{n-2s_{2}}.\]
Consequently, we obtain 
\begin{equation}
\label{PP34}
\| |v(\tau,\cdot)|^{p}\|_{\dot{H}^{s_{1}-1}}\lesssim (1+B_{2}(\tau,0))^{-\frac{n}{2m}p+\frac{n}{4}-\frac{s_{2}-1}{2}}\|(u,v)\|^{p}_{X(t)}.
\end{equation}
Consequently, we get
\begin{eqnarray*}
&&\||D|^{s_{1}-1}u^{nl}_{t}(t,\cdot)\|_{L^{2}(\mathbb{R}^{n})}\\
&\lesssim& \|(u,v)\|_{X(t)}^{p} b_{1}(t)^{-1}(1+B_{1}(t,\tau))^{-\frac{n}{2}\left( \frac{1}{m}-\frac{1}{2}\right)-\frac{s_{1}-1}{2}-1} \\
&&\qquad\times\displaystyle{\int_{0}^{\frac{t}{2}}}b_{1}(\tau)^{-1} (1+B_{1}(\tau,0))^{(-\frac{n}{2m}p+\frac{n}{2m})\beta+\gamma_{1}} d\tau\\
&&+\|(u,v)\|_{X(t)}^{p}b_{1}(t)^{-1} (1+B_{1}(\tau,0))^{(-\frac{n}{2m}p+\frac{n}{4})\beta+\gamma_{1}}\\
&&\qquad\times\displaystyle{\int_{\frac{t}{2}}^{t}} b_{1}(\tau)^{-1}(1+B_{1}(t,\tau))^{-\frac{s_{1}-1}{2}-1}  d\tau\\
&\lesssim& \|(u,v)\|_{X(t)}^{p} b_{1}(t)^{-1}(1+B_{1}(t,\tau))^{-\frac{n}{2}\left( \frac{1}{m}-\frac{1}{2}\right)-\frac{s_{1}-1}{2}-1 },
\end{eqnarray*}
where $\tilde{p}>\frac{2m}{n}\Big(\frac{s_{1}+1+2\gamma_{1}}{2} \Big) +1.$
\section{Appendix} \label{appendix}
Here we state some inequalities which come into play in our proofs.
\begin{Proposition}\label{AP1}
Let $1<p,p_{0},p_{1}<\infty$, $\sigma>0$ and $s\in[0,\sigma).$ Then the following fractional Gagliardo-Nirenberg inequality holds for all $u\in L^{p_{0}} \cap \dot{H} _{p_{1}}^{\sigma}$:
\begin{equation}
\|u\|_{\dot{H} _{p}^{s}}\lesssim\|u\|_{L^{p_{0}}}^{1-\theta}\|u\|_{\dot{H} _{p_{1}}^{\sigma}}^{\theta},
\end{equation}
where \[ \theta=\theta_{s,\sigma}:=\frac{\frac{1}{p_{0}}-\frac{1}{p}+\frac{s}{n}}{\frac{1}{p_{0}}-\frac{1}{p_{1}}+\frac{\sigma}{n}}\,\,\,\mbox{ and}\,\,\,  \frac{s}{\sigma}\leq \theta\leq 1.\]
\end{Proposition}	
For the proof see \cite{hajaieje} and \cite{Christ,Grafakos,Grafakos1,Gulisashvili,Kato,Kenig}.
\begin{Proposition} \label{PropLeibnizrule}
Let us assume $s>0$ and $ 1\leq r\leq\infty, 1<p_{1},p_{2},q_{1},q_{2}\leq \infty$ satisfying the relation
$$	\frac{1}{r}=	\frac{1}{p_{1}}+\frac{1}{p_{2}}=	 \frac{1}{q_{1}}+\frac{1}{q_{2}}.$$
Then the following fractional Leibniz rule holds:
$$	\||D|^{s}(fg)\|_{L^{r}}\lesssim \||D|^{s}f\|_{L^{p_{1}}} \|g\|_{L^{p_{2}}}+\|f\|_{L^{q_{1}}}\||D|^{s}g\|_{L^{q_{2}}},$$
for all $f\in \dot{H} _{p_{1}}^{s}\cap L^{q_{1}}$ and $g\in \dot{H} _{q_{2}}^{s}\cap L^{p_{2}}.$
\end{Proposition}	
For more details concerning fractional Leibniz rule see \cite{Grafakos}.
\begin{Proposition}\label{ap11}
Let us choose $s>0,p>\lceil s \rceil$ and $1<r,r_{1},r_{2}<\infty$ satisfying
$$ \frac{1}{r}=\frac{p-1}{r_{1}}+\frac{1}{r_{2}}.$$
Let us denote by $F(u)$ one of the functions  $|u|^{p},\,\pm|u|^{p-1}u.$
Then the following fractional chain rule holds:
\begin{equation}\label{AP2}
\||D|^{s}F(u)\|_{L^{r}}\lesssim \|u\|_{L^{r_{1}}}^{p-1}	\||D|^{s}u\|_{L^{r_{2}}},
\end{equation}
\end{Proposition}
For the proof see \cite{Palmieri3}.
\begin{Proposition}\label{ap7}
Let $p>1$ and $u\in H^{s}_{m}$, where $s\in(\frac{n}{m},p)$. Then the following estimates hold:
$$	\||u|^{p}\| _{H^{s}_{m}}\lesssim\|u\|_{ H^{s}_{m}}\|u\|^{p-1}_{L^{\infty}},$$
$$	\|u|u|^{p-1}\| _{H^{s}_{m}}\lesssim\|u\|_{ H^{s}_{m}}\|u\|^{p-1}_{L^{\infty}}.$$
\end{Proposition} For the proof see \cite{Sickel}.\\
We can derive from Proposition \ref{ap7} the following corollary.
\begin{Corollary}\label{ap8}
Under the assumptions of Proposition \ref{ap7} it holds
$$	\||u|^{p}\| _{\dot{H}^{s}_{m}}\lesssim\|u\|_{ \dot{H}^{s}_{m}}\|u\|^{p-1}_{L^{\infty}},$$
$$	\|u|u|^{p-1}\| _{\dot{H}^{s}_{m}}\lesssim\|u\|_{ \dot{H}^{s}_{m}}\|u\|^{p-1}_{L^{\infty}}.$$
\end{Corollary}	
For the proof see \cite{Kainane}.
\begin{Lemma}\label{ap15}
Let $0<2s^{*}<n<2s$. Then for any function $f\in\dot{H}^{s^{*}}\cap\dot{H}^{s}$ one has the estimate
$$	\|f\|_{L^{\infty}}\leq\|f\|_{\dot{H}^{s^{*}}}+\|f\|_{\dot{H}^{s}}.$$
\end{Lemma}
For the proof see \cite{D'Abbicco1}.

\end{document}